\documentclass{amsart}
\newtheorem{Lemma}{Lemma}[section]\newcommand{\bel}{\begin{Lemma}}\newcommand{\eel}{\end{Lemma}}
\newtheorem{Proposition}[Lemma]{Proposition}\newcommand{\bprop}{\begin{Proposition}}\newcommand{\eprop}{\end{Proposition}}
\newtheorem{Theorem}[Lemma]{Theorem}\newcommand{\bthe}{\begin{Theorem}}\newcommand{\ethe}{\end{Theorem}}
\newcommand{\bpr}{{\em Proof.\ }}\def\epr{$\bull$\\}
\newtheorem{Remark}[Lemma]{Remark}\newcommand{\beR}{\begin{Remark}\rm}\newcommand{\eeR}{\end{Remark}}
\newtheorem{Definition}[Lemma]{Definition}\newcommand{\bed}{\begin{Definition}}\newcommand{\eed}{\end{Definition}}
\newtheorem{Example}[Lemma]{Example}\newcommand{\bex}{\begin{Example}\rm}\newcommand{\eex}{\end{Example}}
\newtheorem{Corollary}[Lemma]{Corollary}\newcommand{\bcor}{\begin{Corollary}}\newcommand{\ecor}{\end{Corollary}}

\newcommand{\beq}{\begin{equation}}\newcommand{\eeq}{\end{equation}}
\newcommand{\bem}{\begin{displaymath}}\newcommand{\eem}{\end{displaymath}}
\newcommand{\beqa}{\begin{eqnarray}}\newcommand{\eeqa}{\end{eqnarray}}
\newcommand{\bee}{\begin{enumerate}}\newcommand{\eee}{\end{enumerate}}
\newcommand{\bei}{\begin{itemize}}\newcommand{\eei}{\end{itemize}}
\newcommand{\bet}{\begin{tabular}{cccccccc}}\newcommand{\eet}{\end{tabular}}
\newcommand{\bpm}{\begin{pmatrix}}\newcommand{\epm}{\end{pmatrix}}
\newcommand{\bM}{\begin{matrix}}\newcommand{\eM}{\end{matrix}}
\newcommand{\ber}{\begin{array}{l}}\newcommand{\eer}{\end{array}}

\def\bull{\vrule height .9ex width .9ex depth -.1ex }

\def\di{\partial}

\def\isom{\xrightarrow{\sim}}

\def\li{~\\ $\bullet$ }\def\ls{~\\ $\star$ }

\def\cO{\mathcal{O}}

\def\C{\mathbb{C}}

\def\k{k}\def\P{\mathbb{P}}

\def\Z{\mathbb{Z}}

\def\de{\delta}
\def\ep{\epsilon}

  \def\tX{{\tilde{X}}}

\def\suml{\sum\limits}

\def\prodl{\prod\limits}

\def\smin{\setminus}\def\sset{\subset}

\def\?{{\bf ???}}\def\cf{C_{n,r}}
\newcommand{\bin}[2]{\binom{#1}{#2}}
\newcommand{\stir}[2]{\genfrac{\{}{\}}{0pt}{0}{#1}{#2}}

\title{A \MakeLowercase{counterexample to} D\MakeLowercase{urfee's conjecture}}
\author{D\MakeLowercase{mitry} K\MakeLowercase{erner} \MakeLowercase{and} A\MakeLowercase{ndr\'as} N\MakeLowercase{\'emethi}}
\date{\today}

\address{Department of Mathematics, University of Toronto, 40 St. George Street, Toronto, Canada}
\email{dmitry.kerner@gmail.com}

\address{R\'enyi Institute of Mathematics\\Budapest\\Re\'altanoda u. 13\textendash15\\1053\\Hungary}
\email{nemethi@renyi.hu}

\thanks{A.N. is partially supported by OTKA grant  of the
Hungarian Academy of Sciences.}
\thanks{Part of the work was done in Mathematische Forschungsinstitute
Oberwolfach, during D.K.'s stay as an OWL-fellow.}
\keywords{}

\begin{document}\setcounter{secnumdepth}{6} \setcounter{tocdepth}{1}
\begin{abstract}
An old conjecture of Durfee \cite{Durfee1978}
bounds the ratio of two basic invariants of complex isolated
complete intersection surface singularities:
the Milnor number and the singularity (or geometric) genus. We give a counterexample
for the case of non-hypersurface complete intersections, and we
formulate a weaker conjecture valid in arbitrary dimension and codimension.
This weaker bound is asymptotically sharp. In this note we support the validity of
the new proposed inequality by
its verification in certain (homogeneous) cases.
 In our subsequent paper we will prove it for several other cases  and  we
 will provide a more comprehensive  discussion.
\end{abstract}
\maketitle

\section{Introduction}
Let $(X,0)\sset(\C^N,0)$ be the analytic  germ of a complex
 isolated complete intersection singularity (ICIS).
Among various singularity invariants the two most basic ones are:

\vspace{2mm}

\noindent $\bullet$ \  the Milnor number $\mu$, which  measures the  change of the local
topological/homological  Euler characteristic   in smoothing deformation;\\
$\bullet$ \
the singularity genus (called also geometric genus)
 $p_g$, which  measures the change of the local
analytic Euler characteristic (Todd index) when we replace $(X,0)$ by its resolution.

\vspace{2mm}

Their definitions and some of their properties are given in \S\ref{Sec.Background}.

The relation of these two invariants was investigated intensively giving rise to several open problems
  as well. In particular, in \cite{Durfee1978} two  conjectures were formulated:

\vspace{2mm}

\noindent $\bullet$ \
 {\em strong inequality}: if $(X,0)$ is an isolated complete intersection surface singularity, then $\mu\ge 6p_g$;\\
 $\bullet$ \ {\em weak inequality}: if $(X,0)$ is a normal surface singularity (not necessarily ICIS) which admits a smoothing with Milnor number (second Betti number of the fiber) $\mu$,
then $\mu+\mu_0\ge 4p_g$ (where $\mu_0$ is the rank  of the kernel of the intersection form).

\vspace{2mm}

Quite soon a counterexample to the weak conjecture was given in
\cite[page 240]{Wahl 1981} providing a
normal surface singularity (not ICIS) with $\mu=3$, $\mu_0=0$  and $p_g=1$.

On the other hand, the `strong inequality' valid  for an ICIS was believed to be true (though hard to prove) and was verified in many
particular hypersurfaces.  For example,

\vspace{2mm}

\cite{Tomari93} proved $8p_g<\mu$ for $(X,0)$ of multiplicity 2,

\cite{Ashikaga92} proved $6p_g\le\mu-2$ for $(X,0)$ of multiplicity 3,

\cite{Xu-Yau1993} proved $6p_g\le\mu-mult(X,0)+1$ for quasi-homogeneous singularities,

\cite{Nemethi98,Nemethi99} proved $6p_g\le\mu$ for suspension type singularities $\{g(x,y)+z^k=0\}\sset(\C^3,0)$,

\cite{Melle00} proved $6p_g\le\mu$ for absolutely isolated singularities.

\vspace{2mm}

Moreover, for arbitrary dimension, \cite{Yau-Zhang2006} proved the inequality $\mu\ge (n+1)!p_g$
for isolated weighted-homogeneous hypersurface singularities in $(\C^{n+1},0)$. The natural expectation
 was that the same bound $(n+1)!$ holds for any ICIS of any dimension $n$ and any
 codimension.

\vspace{2mm}

In \S\ref{Sec.Counterexample} we give a counterexample to the strong conjecture
in any codimension $r\geq  2$.
In fact, for $r\geq 2$, the  bound  $(n+1)!$ is wrong even asymptotically.
Therefore, in \S\ref{Sec.Weaker.Conjecture} we propose a weaker bound. It is based
on the
Stirling number of the second kind:  $\stir{n+r}{r}:=\frac{1}{r!}\suml^r_{j=0}(-1)^j\bin{r}{j}(r-j)^{n+r}$, cf.
\cite[\S24.1.4]{Abramowitz-Stegun}.

\vspace{2mm}

{\bf Conjecture.} {\em Let $(X,0)\sset(\C^N,0)$ be an ICIS of dimension $n$ and codimension $r=N-n$. Then
\li for $n=2$ and $r=1$ one has  \ $\mu\geq 6p_g$,
\li for $n=2$ and arbitrary $r$ one has  \ $\mu>4p_g$,
\li for $n\ge3$ and fixed $r$ one has  \ $\mu\ge \cf\cdot  p_g$, where the  coefficient $\cf$ is defined by } $$\cf:=\frac{\bin{n+r-1}{n}(n+r)!}{\stir{n+r}{r}r!}.$$

\vspace{2mm}

Note that for any curve singularity (i.e. $n=1$), $p_g$ is the delta invariant $\delta$ and  $\mu\le2 \delta $,
 (with equality exactly for irreducible ones); see example \ref{Ex.Intro.Sing.Genus.For.Curves}.
For $n=2$ the inequality   $\mu(X,0)\ge\cf\cdot  p_g(X,0)$, in general,  is not satisfied
(see  Proposition \ref{Thm.Change.Milnor.Genus.Blowup} and the comments following it).

The bound of the conjecture is asymptotically sharp, i.e. for any fixed $n$ and $r$ there exists a sequence of isolated
complete intersections
 for which the
ratio $\frac{\mu}{p_g}$ tends to $\cf$. In our subsequent longer article \cite{Kerner-Nemethi2011}
we  verify it for several cases;
here we exemplify  only some `homogeneous' situations.
Moreover, we list some elementary properties of the sequence $\{\cf\}_{n,r}$, cf.
\S\ref{Sec.Weaker.Conjecture}. For example,  we show:
\beq
(n+1)!=C_{n,1}\ge C_{n,2}=\frac{(n+2)!(n+1)}{2^{n+2}-2}\ge\cdots\ge \lim_{r\to\infty}\cf=2^n.
\eeq

We wish to thank H. Hamm, A. Khovanskii, M. Leyenson, P. Milman, E. Shustin for advises and
 important discussions.
\section{Background}\label{Sec.Background}
\subsection{The Milnor number}
Let $(X,0)\sset(\C^N,0)$ be an isolated complete intersection of  dimension $n\geq 1$
defined  as the germ of the zero set  of an analytic germ $f:(\C^N,0)\to (\C^r,0)$.
Let $B\sset\C^N$ be
a small ball centered at the origin. By the classical theory (see e.g. \cite{Looijenga-book}) the `Milnor fiber'  $f^{-1}(\ep)\cap B$ ($0<\ep\ll $ the radius of $B$) is smooth and has the homotopy type of a bouquet of $n$--spheres.
Their number $\mu$ is called the Milnor number of $(X,0)$.

For general formulae  of $\mu$ see  e.g. \cite{Milnor-book} and
\cite[\S9.A]{Looijenga-book}.
\subsection{The singularity (geometric)  genus}
Let $(X,0)\sset(\C^N,0)$ be an ICIS as above with $n>1$.
Let $(\tX,E)\to(X,0)$ be one of its resolutions,
i.e. a birational morphism with $\tX$ smooth,
$\tX\smin E\isom X\smin\{0\}$ and $\overline{\tX\smin E}=\tX$.  The singularity
  genus reflects the cohomological non-triviality
of the structure sheaf on the germ $(\tX,E)$, \cite{Artin-66}:
\beq\label{Eq.Def.Singularity.Genus.General}
p_g(X,0):=h^{n-1}\cO_{\tX}.
\eeq
This number does not depend on the choice of the resolution.

\bex\label{Ex.Intro.Sing.Genus.For.Curves}
$\bullet$ If $(X,0)$ is a (reduced) curve singularity then the analogue
 of singularity genus is classical delta invariant:
 $\de=dim_\C \,\cO_{\tX}/\cO_{(X,0)}$. It is related to the Milnor number
 via  $2\de=\mu+\sharp-1$, where $\sharp$ is the
  number of local branches; cf.  \cite{Milnor-book},  \cite{Buchweitz-Greuel-1980},
  \cite{Looijenga1986}.
\li Assume that $(X,0)$ is a normal surface singularity which admits
a smoothing. If $F$ is the Milnor fiber of the smoothing, let $(\mu_+,\mu_0,\mu_-)$
be the Sylvester invariants of the symmetric intersection form in the middle integral
homology  $H_2(F,\Z)$.
Then $2p_g=\mu_0+\mu_+$, \cite[Proposition 3.1]{Durfee1978}.
\li
Let $(X,0)$ be a {\it Gorenstein} normal surface singularity
which admits a smoothing. Then, in fact, the Milnor number of the smoothing  is
independent of the smoothing. Indeed, if $(\tX,E)\to(X,0)$ is any resolution with
relative canonical class $K_{\tX}$,
then $\mu+1=\chi_{top}(E)+K^2_{\tX}+12p_g$.
This formula was proved  in \cite{Laufer1977} for hypersurface surface singularities
  and in \cite{Looijenga1986} for normal Gorenstein smoothable singularities.
\li
Let $(X,0)\sset(\C^{n+1},0)$ be an isolated hypersurface singularity. Let $Sp(X,0)\sset(-1,n)$ be its singularity
spectrum. Then  $p_g$ is the cardinality of $Sp\cap(-1,0]$, cf. \cite{Saito1981,Steenbrink-1983}.
\eex

For the general introduction to singularities we refer to \cite{AGLV}, \cite{Dimca92}, \cite{Looijenga-book},
\cite{Seade}.
\section{The counterexample}\label{Sec.Counterexample}
Let $Y=\{f_1=\cdots=f_r=0\}\sset\P^{N-1}$ be a smooth
projective complete intersection defined by homogeneous
polynomials of degrees $deg(f_i)=p_i$. Define the corresponding ICIS as the cone over $Y$: $(X,0)=Cone(Y)\sset(\C^N,0)$.
\bprop\label{Thm.Milnor.Genus.For.ICIS}
(1) \ \ $\mu=\Big(\prodl^r_{i=1}p_i\Big)\suml^n_{j=0}(-1)^j
\bigg(\suml_{\substack{(k_1,..,k_r),\ k_i\ge0\\\suml_i k_i=n-j}}
\prod_{i=1}^r(p_i-1)^{k_i}\bigg)-(-1)^n$,
\\(2) \ \ \ \ \ \
$p_g=\suml_{\substack{(k_1,..,k_r),\ k_i\ge0\\\suml_i k_i=n}}\prodl^r_{i=1}\bin{p_i}{k_i+1}$.
\eprop
For a systematic study of the Milnor number of weighted homogeneous complete intersections
the reader is invited to consult
 \cite{Greuel1975,Greuel-Hamm1978}. Additionally, several other formulae can be found in the literature, see e.g.
\cite{Hamm1986,Hamm2011}. For formulae regarding  the geometric genus we refer to
\cite{Khovanskii1978,Morales1985,Hamm2011}. These  formulae usually are rather different
 than ours considered above. Nevertheless, the above
expressions can be derived from them: here for the Milnor number we will use \cite{Greuel-Hamm1978}, and for the
geometric genus \cite{Morales1985}.

\vspace{2mm}

\bpr (1) We will determine the Euler characteristic $\chi=(-1)^n\mu +1$ of the Milnor fiber.
For a power series $Z:=\sum_{i\geq 0}a_ix^i$ write $[Z]_n$ for the coefficient $a_n$ of $x^n$. Also, set
$P:=\prod_{i=1}^rp_i$. Then, by formula 3.7(c) of  \cite{Greuel-Hamm1978}
$$\chi=P\cdot \Big[\frac{(1+x)^N}{\prod_i(1+p_ix)}\Big]_n.$$
Rewrite $1+p_ix$ as $(1+x)(1-\frac{(1-p_i)x}{1+x})$, hence
$$\Big[\frac{(1+x)^N}{\prod_i(1+p_ix)}\Big]_n=
\Big[(1+x)^n\cdot \prod_i\sum_{k_i\geq 0}\Big(\frac{(1-p_i)x}{1+x}\Big)^{k_i}\Big]_n=$$
$$\Big[ \sum _{k_1\geq 0,..,k_r\geq 0}x^{\sum k_i}(1+x)^{n-\sum k_i}\prod_i(1-p_i)^{k_i}\Big]_n=
 \sum _{k_1\geq 0,..,k_r\geq 0\atop \sum k_i\leq n}\prod_i(1-p_i)^{k_i}.
$$
(2) By Theorem 2.4 of \cite{Morales1985} (and  computation of the lattice point under the
`homogeneous Newton diagram') we get
\beq\label{pg}
p_g=\bin{\sum_kp_k}{N}-\sum_{1\leq i\leq r}\bin{(\sum_kp_k)-p_i}{N}+
\sum_{1\leq i<j\leq r}\bin{(\sum_kp_k)-p_{i}-p_{j}}{N}-\ldots
\eeq
Using the  Taylor expansion $\frac{1}{(1-z)^{N+1}}=\sum_{l\geq N}\bin{l}{N}z^{l-N}$,
the right hand side of (\ref{pg})  is $\Big[\frac{\prod_i (1-z^{p_i})}{(1-z)^{N+1}}\Big]_{\sum_kp_k-N}$.
Thus
$$
p_g=rez_{z=0}F(z), \ \ \mbox{where} \ \
F(z):=\frac{\prod_i (1-z^{p_i})}{z^{\sum p_i-N+1}(1-z)^{N+1}}.
$$
Note that $rez_{z=\infty}F(z)=0$ since  $F(1/z)/z^2$ is regular at zero. Since
$F(z)$ has poles at $z=0$ and $z=1$ only, and  $\sum_{pt} rez_{z=pt}F(z)=0$,
we have $p_g=-rez_{z=1}F(z)$.
By the change of variables $z\mapsto 1/z$ we get
$$
p_g=rez_{z=1}\frac{\prod_{i=1}^r (z^{p_i}-1)}{(z-1)^{N+1}}.
$$
Since $z^p-1=\sum_{k\geq 1}\binom{p}{k}(z-1)^k$, we obtain
$$
p_g
=rez_{z=1}\frac{1}{(z-1)^{n+1}}\cdot \prod_{i=1}^r \sum_{k_i\geq 0}\binom{p_i}{k_i+1}(z-1)^{k_i}=
\sum_{k_1,\ldots, k_r\geq 0\atop k_1+\cdots +k_r=n}\binom{p_1}{k_1+1}\cdots \binom {p_r}{k_r+1}.
$$ This ends the proof. \epr

\bex
Consider the particular case $p_1=\cdots=p_r=p$. Then the formulae of Proposition \ref{Thm.Milnor.Genus.For.ICIS}
read as
\beq\label{eq:pp}\ber
\mu=(-1)^n\Big(p^r\sum_{j=0}^n(1-p)^j\bin{j+r-1}{j}-1\Big), \ \mbox{(see also \cite[3.10(b)]{Greuel-Hamm1978})},
\\\\p_g=\sum_{\substack{(k_1,..,k_r),\ k_i\ge0\\\suml_i k_i=n}}\prod^r_{i=1}\bin{p}{k_i+1}.
\eer\eeq
In some low dimensional cases we have:

\ls For $n=1$: $p_g=\frac{r(p-1)p^r}{1!2^1}$
\ls For $n=2$: $p_g=\frac{r(p-1)p^r}{2! 2^2}\Big(r(p-1)+\frac{p-5}{3}\Big)$
\ls For $n=3$: $p_g=\frac{r(p-1)p^r}{3!2^3}(p r-2 - r) (p r-3 + p - r)$.

\vspace{2mm}

Analyzing (\ref{eq:pp}) one sees that in the case of $\mu$, the leading term in $p$ (for $p$ large)
 comes from the last summand ($j=n$) and  $\mu= p^N\bin{N-1}{n}+O(p^{N-1})$.

The leading term for $p_g$ is more complicated. For $p\gg0$ one has:
$$p_g=p^N\cdot\suml_{\substack{(k_1,..,k_r),\ k_i\ge0\\\suml_i k_i=n}}\prodl^r_{i=1}\frac{1}{(k_i+1)!}+O(p^{N-1}).$$

Thus, asymptotically,
$\frac{\mu}{p_g}=\frac{\bin{N-1}{n}}{\suml_{\substack{(k_1,..,k_r),\ k_i\ge0\\\suml_i k_i=n}}
\prodl^r_{i=1}\frac{1}{(k_i+1)!}}+O(\frac{1}{p})$, which in low dimensions gives:
\beq
n=N-r=2:\quad \frac{\mu}{4p_g}=\frac{r+1}{r+\frac{1}{3}}+O(\frac{1}{p});\quad\quad\quad
n=N-r=3:\quad \frac{\mu}{8p_g}=\frac{r+2}{r}+O(\frac{1}{p}).
\eeq
Since $4(r+1)/(r+\frac{1}{3})< 6$  for all $r\geq 2$, in the case $n=2$ the
strong Durfee bound is violated even asymptotically --- that is, for any $p$ sufficiently large ---
whenever $r\geq 2$. Similarly,
for $n=3$, $8(r+2)/r<4! $ whenever $r\geq 2$.
\eex
\section{The new conjecture with the new weaker bound}\label{Sec.Weaker.Conjecture}
In the previous section we have seen that
 for the singularity which is the cone over a smooth projective
complete intersection with $p_1=\cdots=p_r$, the ratio $\mu/p_g$ tends
to the numerical factor
\beq
\cf:=\frac{\bin{N-1}{n}}{\suml_{\substack{(k_1,..,k_r),\ k_i\ge0\\\suml_i k_i=n}}\prodl^r_{i=1}\frac{1}{(k_i+1)!}}.
\eeq
Next we list some properties of these numbers.
\bel
(I) \ $\cf=\frac{\bin{n+r-1}{n}(n+r)!}{\stir{n+r}{r}r!}$

\vspace{2mm}

(II) \ $(n+1)!=C_{n,1}\ge C_{n,2}\ge\cdots\ge \lim_{r\to\infty}\cf=2^n$.
\eel
\bpr
(I) \ Let $S_{n,r}:=\sum_{\substack{\{k_i\geq 0\}\\\sum_i k_i=n}}\prod^r_{i=1}\frac{1}{(k_i+1)!}$.
First we prove
\beq\label{eq:new2}
\frac{r^{n+r}}{(n+r)!}=S_{n,r}+rS_{n+1,r-1}+\bin{r}{2}S_{n+2,r-2}+\cdots.\eeq
Consider the expansion
$r^{n+r}=(1+\cdots+1)^{n+r}=\sum_{\substack{k_1,..,k_r\ge0}}(\prod^r_{i=1}1^{k_i})\frac{(n+r)!}{\prod^r_{i=1}k_i!}$.
This gives
\beq\label{eq:new}
\frac{r^{n+r}}{(n+r)!}=\sum_{\substack{k_1,..,k_r\ge0\\ \sum_ik_i=n+r}}
\frac{1}{\prod^r_{i=1}k_i!}.
\eeq
For each $0\leq j\leq r$ and subset $I(j)\in\{1,\cdots,r\}$ of cardinality $j$  set $K_{I(j)}:=\{\underline{k}=(k_1,\ldots, k_r)\ \mbox{with} \
 k_i=0 \ \mbox{for} \ j\in I(j)\}$.
Then  the sum $\sum_{\underline{k}}$  from the right hand side of (\ref{eq:new})
can be replaced by
 $\sum_{j=0}^r \sum_{I(j)}\,\sum_{\underline{k}\in K_{I(j)}}$.
 Since for each fixed $I(j)$
 $$\sum_{\underline{k}\in K_{I(j)}}\frac{1}{\prod^r_{i=1}k_i!}=S_{n+j,r-j}$$
 (by a change $k_j\mapsto k_j-1$ for those indices when $k_j\not=0$),
 and for each $j$ there are $\binom{r}{j}$ subsets $I(j)$  of cardinality $j$, (\ref{eq:new2})  follows.
\\

Now, we prove that $S_{n,r}=\frac{\sum^r_{i=0}(r-i)^{n+r}(-1)^i\bin{r}{i}}{(n+r)!}$.
 Note that the recursion (\ref{eq:new2})
 defines $S_{n,r}$ uniquely. Hence it is enough to verify
that the proposed solution satisfies this recursive relation; that is,
$r^{n+r}=\sum^r_{j=0}\bin{r}{r-j}\sum^j_{i=0}(j-i)^{n+r}(-1)^i\bin{j}{i}$.
Set  $\tilde{j}=j-i$. Then the sum can be rewritten as
$\sum^r_{\tilde{j}=0}\tilde{j}^{n+r}\bin{r}{\tilde{j}}\sum^{r-\tilde{j}}_{i=0}(-1)^i\bin{r-\tilde{j}}{i}$.
Note that $\sum^{r-\tilde{j}}_{i=0}(-1)^i\bin{r-\tilde{j}}{i}=(1-1)^{r-\tilde{j}}=0$ unless $r=\tilde{j}$. Hence
the above formula for $S_{n,r}$ follows.

This formula can be compared with that one satisfied by the Stirling numbers (see introduction), hence
we get $S_{n,r}=\stir{n+r}{r}\frac{r!}{(n+r)!}$, which proves  the statement for $\cf$.
\\(II) \
The limit of $\cf$ is computed using the asymptotics of Stirling numbers of the second kind,
\cite[\S24.1.4]{Abramowitz-Stegun}:
$\stir{n+r}{r}\sim\frac{r^{2n}}{2^nn!}$. This gives: $\cf\sim 2^n\frac{(n+r-1)!(n+r)!}{(r-1)!r!r^{2n}}\to 2^n$.

The following proof showing that   $\cf$ is non-increasing  sequence in $r$ was communicated to us
by D. Moews.
Write the inequality $\cf\ge C_{n,r+1}$ in terms of Stirling numbers of the second
kind:
\beq
\stir{n+r+1}{r+1}r(r+1)\ge\stir{n+r}{r}(n+r)(n+r+1).
\eeq
Set $N:=n+r$ as usual, and write the inequality in terms of generating functions:
\beq
\frac{r(r+1)}{x}\sum_{N\ge r}\stir{N+1}{r+1}\frac{x^{N+1}}{(N+1)!}\succeq x\di_x\sum_{N\ge r}\stir{N}{r}\frac{x^N}{N!}
\eeq
with the convention: $\sum a_n x^n\succeq\sum b_n x^n$ if and only if  $ a_n\ge b_n$ for any $n\geq 0$.

Recall that Stirling numbers are the coefficients of Taylor series:
 $\frac{(e^x-1)^r}{r!}=\sum_{N\ge r}\stir{N}{r}\frac{x^N}{N!}$. Therefore the inequality to be proved can be
 rewritten as:
$\frac{(e^x-1)^{r+1}}{(r-1)!x}\succeq  \frac{x e^x(e^x-1)^{r-1}}{(r-1)!}$.
We claim that this inequality will follow from $(e^x-1)^2\succeq  x^2e^x$. Indeed, if $\sum a_nx^n\succeq\sum b_n x^n$
and the Taylor expansion of $g(x)$ has only positive coefficients, then $g(x)\sum a_nx^n\succeq g(x)\sum b_n x^n$.
Hence it is enough to prove $(e^x-1)^2\succeq  x^2e^x$. For the same reason, this last statement
 will follow from $(e^x-1)\succeq  x e^{\frac{x}{2}}$,  which is immediate.
\epr

The coefficients $C_{n,r}$ satisfies the following identities and inequalities (some more will be listed in
\cite{Kerner-Nemethi2011}, where  the general inequality  $\mu\ge \cf p_g$  will also be discussed).

\bprop\label{Thm.Change.Milnor.Genus.Blowup}
Let $(X,0)=\cap^r_{i=1}\{f_i(x_1,\ldots,x_N)=0\}\sset(\C^N,0)$ be an isolated complete intersection
singularity, where each $f_i$ is homogeneous  of degrees $p_i$. Set $P:=\prod^r_{i=1}p_i$.
Then:

\vspace{2mm}

\noindent For $n=1$: $\mu+P-1=C_{1,r}\cdot p_g=2p_g$. (Here $p_g:=\delta$.)
\\For $n=2$:
$\mu+P\cdot \Big(\frac{r-1}{3r+1}\sum^r_{i=1}(p_i-1)-\sum_{i<j}\frac{(p_i-p_j)^2}{3r+1}-1\Big)+1=C_{2,r}\cdot p_g=4\frac{r+1}{r+\frac{1}{3}}\cdot p_g$.
\\For $r=1$: $\mu-C_{n,1}\cdot p_g=(p-1)^N-\frac{p!}{(p-N)!}\ge0$. In particular, $\mu\geq(n+1)!\cdot p_g$.
\\For $n=3$: $\mu>C_{3,r}\cdot p_g$.
\eprop
We wish to emphasize the equation $\mu+P\cdot E+1=C_{2,r}\cdot p_g$ in the case $n=2$. Note that
the expression $E$ is $-1$ (hence negative) if $r=1$, but for $r\geq 2$ and some choices of $p_i$'s
(e.g. whenever they are all equal) $E$ might be positive.

\vspace{2mm}

\bpr By Proposition \ref{Thm.Milnor.Genus.For.ICIS} we have:
\beq
\mu=P\cdot \sum^n_{j=0}(-1)^j\bigg(\sum_{\substack{\{k_i\ge0\},\\\sum_i k_i=n-j}}
\prod_{j=1}^r(p_i-1)^{k_i}\bigg)-(-1)^n,\quad\quad
p_g=P\cdot \sum_{\substack{\{k_i\ge0\}\\\sum_i k_i=n}}\prod^r_{i=1}\bin{p_i-1}{k_i}\frac{1}{k_i+1}.
\eeq
We can assume $p_i>1$, otherwise some of defining hypersurfaces are hyperplanes, so
 this reduces to a singularity of {\em the same dimension} and smaller codimension.
\\{\bf Case n=1.} Use
\beq
\mu=P\cdot\Big(\sum^r_{i=1}(p_i-1)-1\Big)+1,\quad\quad
p_g=P\cdot \sum^r_{i=1}\frac{(p_i-1)}{2},\quad\quad C_{1,r}=2.
\eeq
\\{\bf Case n=2.} We have
\beq\ber
\mu=P\cdot
\Big(\suml_{i}(p_i-1)(p_i-2)+\suml_{i<j}(p_{i}-1)(p_{j}-1)+1\Big)-1,\\
p_g=P\cdot
\Big(\suml_{i}\bin{p_i-1}{2}\frac{1}{3}+\suml_{i<j}\frac{(p_{i}-1)(p_{j}-1)}{4}\Big),
\quad\quad C_{2,r}=4\frac{r+1}{r+\frac{1}{3}}.
\eer\eeq
Hence
\beq\ber
\mu-C_{2,r}p_g=P\cdot \Bigg(\suml_{i}(p_i-1)(p_i-2)\frac{(r-1)}{3r+1}-
\frac{2}{3r+1}\suml_{i<j}(p_{i}-1)(p_{j}-1)+1\Bigg)-1=\\
=P\cdot \Bigg(
\suml_{i<j}\frac{(p_{i}-1)^2-2(p_{i}-1)(p_{j}-1)+(p_{j}-1)^2}{3r+1}-\suml_{i}(p_i-1)\frac{(r-1)}{3r+1}+1\Bigg)-1,
\eer\eeq
giving the needed equality. Note that from here one also obtains for $r=1$ and $n=2$:
\beq
6p_g=\mu-P+1\leq \mu.
\eeq
\\{\bf Case $r=1$} is obvious.
\\{\bf Case $n=3$.}  It is enough to  prove:
\beq
\sum^3_{j=0}(-1)^j\bigg(\sum_{\substack{\{k_i\ge0\},\\\sum_i k_i=3-j}}\prod_{j=1}^r(p_i-1)^{k_i}\bigg)>
\frac{8(r+2)}{r}\sum_{\substack{\{k_i\ge0\}\\\sum_i k_i=3}}\prod^r_{i=1}\bin{p_i-1}{k_i}\frac{1}{k_i+1}
\eeq
Present the difference as follows
\beq\ber
\underbrace{\sum_{\substack{\{k_i\ge0\},\\\sum_i k_i=3}}\prod_{j=1}^r(p_i-1)^{k_i}\Big(1-\frac{8(r+2)}{r(k_i+1)!}\Big)}_{Part\ I}
+\underbrace{\frac{8(r+2)}{r}\sum_{\substack{\{k_i\ge0\}\\\sum_i k_i=3}}\prod^r_{i=1}
\Big(\frac{(p_i-1)^{k_i}}{(k_i+1)!}-\bin{p_i-1}{k_i}\frac{1}{k_i+1}\Big)}_{Part\ II}+\\
+\underbrace{\sum^3_{j=1}(-1)^j\bigg(\sum_{\substack{\{k_i\ge0\},\\\sum_i k_i=3-j}}\prod_{j=1}^r(p_i-1)^{k_i}\bigg)}_{Part\ III}.
\eer\eeq
First we prove that $Part\ I\ge0$.
\beq
Part\ I=\frac{2r-2}{3r}\sum_i(p_i-1)^3+\frac{r-4}{3r}\sum_{i\neq j}(p_{i}-1)^2(p_{j}-1)-
\frac{2}{r}\sum_{i<j<k}(p_{i}-1)(p_{j}-1)(p_{k}-1).
\eeq
Now use the elementary inequalities (the algebraic mean compared to the geometric mean)
\beq\ber
\sum a^n_i=\suml_{i\neq j}\frac{l a^n_{i}+(n-l)a^n_{j}}{n(r-1)}\ge\frac{1}{r-1}\suml_{i\neq j}a^l_{i}a^{n-l}_{j},
\quad
\suml_{\substack{i\neq j\\2k<n}}a^k_{i}a^{n-k}_{j}=\suml_{i\neq j\neq m}
\frac{(n-k-1)a^k_{i}a^{n-k}_{j}+a^k_{i}a^{n-k}_{m}}{(n-k)(r-2)}
\ge\frac{\suml_{i\neq j\neq m}a^k_{i}a^{n-k-1}_{j}a_{m}}{r-2}\\
\eer\eeq
to get:
\beq
Part\ I\ge \frac{2+r-4}{3r}\sum_{i\neq j}(p_{i}-1)^2(p_{j}-1)-
\frac{2}{r}\sum_{i<j<k}(p_{i}-1)(p_{j}-1)(p_{k}-1)\ge0.
\eeq
Now simplify $Part\ II$:
\beq
Part\ II=\frac{8(r+2)}{r}\Big(\sum_i\frac{(p_i-1)(3p_i-5)}{4!}+\sum_{i\neq j}
\frac{(p_{i}-1)}{3!}\frac{(p_{j}-1)}{2}\Big).
\eeq
Thus, by direct check one gets: $Part\ II-Part\ III>0$. Hence the statement.
\epr

\vspace{2mm}

Finally, we add a slightly weaker inequality, but which is valid for any $n$ (including $n=2$) and any $r$.
Since $\bin{N-1}{n}$ is the number of terms in the sum from the denominator of $C_{n,r}$,
we obtain that $C_{n,r}\geq
\underset{\substack{\sum k_i=n\\k_i\ge0}}{\min}\Big\{\prod^r_{i=1}(k_i+1)!\Big\}$.
\bprop\label{Thm.Bound.mu.pg.ICIS.absolutely.isolated}
Let $(X,0)$ as in Lemma \ref{Thm.Change.Milnor.Genus.Blowup} and $n\ge2$.
Then
$$\mu\ge\underset{\substack{\sum k_i=n\\k_i\ge0}}{\min}\Big\{\prod^r_{i=1}(k_i+1)!\Big\}\cdot p_g\geq 2^n\cdot p_g.$$
Moreover, if
$n>r$, let $n=n_1 r+r_1$ for $n_1,r_1\ge0$, $r_1<r$.
Then $\mu>\Big((n_1+1)!\Big)^{r-r_1}\Big((n_1+2)!\Big)^{r_1}p_g$.

In particular, for $r=1$:  $\mu>(n+1)! p_g$. For $r=2$ and $n$-odd: $\mu\ge (\lfloor\frac{n}{2}\rfloor+1)!(\lfloor\frac{n}{2}\rfloor+2)!p_g$.
For $r=2$ and $n$-even:  $\mu\ge \Big((\lfloor\frac{n}{2}\rfloor+1)!\Big)^2p_g$.
\eprop\
\\\bpr {\bf Step 1.}  We claim that (by
Proposition \ref{Thm.Milnor.Genus.For.ICIS}) it is enough to prove the inequality
\beq
L.H.S.:=P\cdot \bigg(\sum_{j=0}^n(-1)^j\sum_{\substack{\{k_i\ge0\}\\\sum_j k_j\le n-j}}
\prod_{i=1}^r(p_i-1)^{k_i}\bigg)-(-1)^n>
\sum_{\substack{k_1,\ldots, k_r\geq 0\\k_1+\cdots +k_r=n}}\Big(\prod^r_{i=1}\bin{p_i}{k_i+1}(k_i+1)!\Big)=:R.H.S.
\eeq
Clearly this imply the first inequality of the proposition, while the second
follows from $(k+1)!\geq 2^k$.

{\bf Step 2.} The left hand side can be written as
\beq
L.H.S.=(p_r-1)(\prodl^{r-1}_{i=1}p_i)C_n(p_1,..,p_{r})+
(p_{r-1}-1)(\prodl^{r-2}_{i=1}p_i)C_n(p_1,..,p_{r-1})+..+(p_1-1)(p_1-1)^n,
\eeq
where $C_n(p_1,..,p_{r}):=\suml_{\substack{\{k_i\ge0\}\\\sum_i k_i=n}}\Big(\prodl^r_{i=1}(p_i-1)^{k_i}\Big)$.
This is expansion in terms with decreasing $r$, hence it is natural to expand the right hand side similarly.

Let $D_n(p_1,..,p_r):=\suml_{\substack{k_1,\ldots, k_r\geq 0\\k_1+\cdots +k_r=n}}\prod^r_{i=1}\bin{p_i-1}{k_i}(k_i)!$.
So $D_1(p_1)=\bin{p_1-1}{n}n!$. Write the $R.H.S.$ as follows:
\beq
\Big((\prod^r_{i=1}p_i)D_n(p_1,..,p_r)-(\prod^{r-1}_{i=1}p_i)D_n(p_1,..,p_{r-1})\Big)+
\Big((\prod^{r-1}_{i=1}p_i)D_n(p_1,..,p_{r-1})-(\prod^{r-2}_{i=1}p_i)D_n(p_1,..,p_{r-2})\Big)+...+p_1 D_1(p_1).
\eeq
Thus it is enough to prove the inequality for each pair of terms in these expansions:
\beq\label{eq:24}
\forall r\ge1,\ n\ge2:\quad (p_r-1)(\prodl^{r-1}_{i=1}p_i)C_n(p_1,..,p_{r})\ge
(\prod^r_{i=1}p_i)D_n(p_1,..,p_r)-(\prod^{r-1}_{i=1}p_i)D_n(p_1,..,p_{r-1}).
\eeq
For example, for $r=1$ we have: $(p_1-1)^{n+1}> p_1\bin{p_1-1}{n}n!$.
Present the right hand side of (\ref{eq:24}) in the form
$$(p_r-1)(\prod^{r-1}_{i=1}p_i)D_n(p_1,..,p_r)+(\prod^{r-1}_{i=1}p_i)\Big(D_n(p_1,..,p_r)-D_n(p_1,..,p_{r-1})\Big).$$
Note that
\beq
D_n(p_1,..,p_r)-D_n(p_1,..,p_{r-1})=
\suml_{\substack{\{k_i\ge0\},\  k_r>0\\k_1+\cdots +k_r=n}}\prod^r_{i=1}\bin{p_i-1}{k_i}(k_i)!=
(p_r-1)D_{n-1}\big(p_1,..,p_{r-1},(p_r-1)\big).
\eeq
Hence we have to prove the inequality:
 $C_n(p_1,..,p_{r})\ge D_n(p_1,..,p_r)+D_{n-1}\big(p_1,..,p_{r-1},(p_r-1)\big)$.
In fact, we will prove by induction on $r$ the stronger inequality:
\beq
C_n(p_1,..,p_{r})\ge D_n(p_1,..,p_r)+D_{n-1}\big(p_1,..,p_r\big).
\eeq
{\bf Step 3.} Both parts are summations over $(k_1,..,k_r)$, expand them in $k_r$. Then we have to prove:
\beq
\sum^n_{k_r=0}\Bigg((p_r-1)^{k_r}C_{n-k_r}(p_1,..,p_{r-1})-\bin{p_r-1}{k_r}k_r!D_{n-k_r}(p_1,..,p_{r-1})\Bigg)\ge
\sum^n_{k_r=1}\bin{p_r-1}{k_r-1}(k_r-1)!D_{n-k_r}(p_1,..,p_{r-1}).
\eeq
For $k_r\ge2$ one has: $(p_r-1)^{k_r}\ge\bin{p_r-1}{k_r}k_r!+\bin{p_r-1}{k_r-1}(k_r-1)!$, by direct check.
 And $C_j(p_1,..,p_{r-1})\ge D_j(p_1,..,p_{r-1})$, for all values of $j,r$. Therefore we only need to check
 the terms for $k_r=0,1$, i.e. it is enough to prove:
\beq
(p_r-1)\Bigg(C_{n-1}(p_1,..,p_{r-1})-D_{n-1}(p_1,..,p_{r-1})\Bigg)+C_{n}(p_1,..,p_{r-1})-D_n(p_1,..,p_{r-1})\ge
D_{n-1}(p_1,..,p_{r-1}).
\eeq
Again, $C_{n-1}(p_1,..,p_{r-1})\ge D_{n-1}(p_1,..,p_{r-1})$, so the initial inequality reduces to
$C_{n}(p_1,..,p_{r-1})\ge D_n(p_1,..,p_{r-1})+D_{n-1}(p_1,..,p_{r-1})$. This completes the induction step from
$(r-1)$ to $r$.

Finally, for $r=1$ the initial inequality is: $(p_1-1)^n\ge\bin{p_1-1}{n}n!+\bin{p_1-1}{n-1}(n-1)!$. By direct check
it is true for $n\ge2$.
\epr

\end{document}